\documentclass[12pt]{article}
\usepackage{amsfonts}
\usepackage{amsmath}
\usepackage{amscd}
\usepackage{latexsym}

\newtheorem{thm}{{\sc Theorem}}

\newcommand{\qqed}{\hspace*{\fill} $\Box$}

\title{Simultaneous minimal model of homogeneous toric deformation}
\author{Daisuke Matsushita
        }
\date{}

\begin{document}
\maketitle
\begin{abstract}
 For a flat
 family of Du Val singularities, we can take a simultaneous
 resolution after finite base change. It is an interesting 
 problem to consider this analogy for a flat 
 family of higher dimensional
 canonical singularities. In this note, 
 we consider an existence of
 simultaneous terminalization for K.~Altmann's 
 homogeneous toric deformation whose central fibre is an
 affine Gorenstein toric singularity.
 We obtain examples that
 there are no simultaneous terminalization even after finite base
 change and give a sufficient condition for an existence of
 simultaneous terminalization.
 Some examples of 4-dimensional flops are obtained as an
 application.
\end{abstract}
\section{Introduction}
 For a flat family of surfaces $f : X  \to S$, 
 it is an interesting problem exists
 a birational morphism $\tau : \widetilde{X} \to X$ such that
\begin{enumerate}
 \item  $f \circ \tau$ is a flat morphism
 \item  The fibre $\widetilde{X}_s := (f\circ \tau)^{-1}(s)$
        $(s \in S)$ is the minimal resolution of 
        $X_s$.
\end{enumerate}
 If there exists such $\tau$, 
 we say  $f$ admits  a simultaneous resolution.
 Let $f : X \to S$ be a flat morphism whose fibres have
 only Du Val singularities. Then
 Brieskorn \cite{brie1,brie2} and Tyurina\cite{tyur} show
 that 
 there exists a finite surjective morphism $S' \to S$ and
 a flat morphism
 $f' : X\times_{S}S' \to S'$ admits 
 a simultaneous resolution.
 It is a key fact for Minimal Model Theory. 
 Thus it is natural to consider an analogy for
 higher dimensional singularity. 
 According to Minimal Model theory, 
 it is suitable to consider an existence of
 ``simultaneous terminalization'' for
 a flat family of higher dimensional singularity.

\noindent
{\sc Definition 1. \quad}
 For a flat morphism $f : X \to S$ whose
 fibres are $n$-folds $(n \ge 3)$, we say $f$ admits a 
 simultaneous terminalization if there exists a birational
 morphism $\tau : \widetilde{X} \to X$ which satisfies the
 following conditions:
\begin{enumerate}
 \item $f \circ \tau$ is a flat morphism.
 \item The fibre $\widetilde{X}_s := (f\circ \tau)^{-1}(s)$ $(s \in S)$
       has only terminal singularity.
 \item $K_{\widetilde{X}_s}$ is $\tau$-nef.
\end{enumerate}
 Recently 
 K.~Altmann constructed in
 \cite[Definition 3.1]{homogeneous}
 affine toric varieties which are called
 ``homogeneous toric deformation''.
 This toric varieties 
  can be described 
 many flat families of toric singularities such as
 versal deformation space of an isolated Gorenstein toric singularity
 \cite[Theorem 8.1]{alt}.
 In this note, using toric Minimal Model Theory, 
 we investigate the existence of simultaneous terminalization of
 Gorenstein homogeneous toric deformation which is a 
 homogeneous toric deformation whose  central fibre is 
 an affine Gorenstein toric variety. 
 Our main results are as follows:
\begin{thm}
 There exists a 
 Gorenstein homogeneous toric deformation whose fibre 
 is $n$-dimensional $(n \ge 3)$ and which 
 admits no simultaneous
 terminalization even after finite base change.
\end{thm}
 We consider the condition of an existence of simultaneous 
 terminalization, and obtain the following results:
\begin{thm}
 A Gorenstein homogeneous toric deformation $f : X \to {\mathbb C}^{m}$ 
 admits a simultaneous terminalization 
 if $X$ has a crepant desingularization.
\end{thm}
 We recall 
 the definition of 
 homogeneous toric deformation 
 by K.~Altmann, and Gorenstein homogeneous
 toric  deformation in section 2.
 Theorems 1 and 2 will be proved in Sections 3 and 4, respectively.
 Some examples of 4-dimensional flops 
 are obtained as an application.

\section{Homogeneous toric deformation}
\noindent
 First we introduce the definition homogeneous
 toric deformation by K.~Altmann.

\noindent
{\sc Definition 2. \quad} A flat morphism  $f : X \to {\mathbb C}^{m}$ 
is called a homogeneous toric deformation if the following conditions are 
satisfied:
\begin{enumerate}
 \item $X := {\rm Spec}{\mathbb C}[\sigma^{\vee}\cap M]$ 
       is an affine toric variety.
 \item $f$ is defined by $m$ equations such that
       $x^{r_i} - x^{r_0}$ where $r_i \in \sigma^{\vee} \cap M$, 
       $(0 \le i \le m)$.
 \item Let $L := \oplus^{i = m}_{i = 1}{\mathbb Z}(r_{i} - r_{0})$ 
       be the sublattice of $M$.
       The central fibre $Y := f^{-1}(0,\cdots , 0)$ is isomorphic to
       ${\rm Spec}{\mathbb C}[\bar{\sigma}^{\vee}\cap \bar{M}]$ where
       $\bar{\sigma} = \sigma \cap L^{\perp}$ and $\bar{M} := M/L$.
 \item $i : Y \to X$ sends the
       closed orbit in $Y$  isomorphically onto the closed orbit in
       $X$. 
\end{enumerate}
 In this note, we consider a homogeneous toric deformation
 with some additional conditions:

\noindent
{\sc Definition 3. \quad}
 We call homogeneous toric deformation $f: X \to {\mathbb C}^{m}$ a
 Gorenstein homogeneous toric deformation if it 
 satisfies the following two conditions:
\begin{enumerate}
 \item $Y$ is an affine Gorenstein singularity.
 \item The restriction $-r_i$, $(0 \le i \le m)$ to $\bar{M}$ 
       determines $K_Y$.
\end{enumerate}
{\sc Examples.}
 \begin{enumerate}
  \item Most plain example is $f : {\mathbb C}^{2}(x,y) \to {\mathbb C}(t)$.
        In this case, $f$ is defined by $x - y = t$.
  \item Let $g : {\cal X} \to S$ be a versal deformation space of
        $A_n$ singularity. It can be described as follows:
\[
   {\cal X} = (xy + z^{n+1} + t_1 z^{n-1} + \cdots + t_{n-1}z + t_n = 0)
    \subset {\mathbb C}^{n+3} \to {\mathbb C}^{n}(t_1 , \cdots t_n).
\]
        Let $\alpha_{i}$, $(0 \le i \le n)$ be elementary symmetric
        polynomials on ${\mathbb C}^{n+1}(s_0 , \cdots  s_n)$ and
        $H$ a hyperplane such that $\sum_{i=0}^{n}s_i = 0$.
        We take a  base change by 
        $\alpha_i (s_0 , \cdots , s_{n})=t_i  :
        H \to {\mathbb C}^{n}$, 
\[
\begin{CD}
 {\cal X}\times_{{\mathbb C}^{n}} H @>>> {\cal X} \\
 @V{f}VV                  @VV{g}V\\
 H                @>>> {\mathbb C}^{n}.
\end{CD}
\]
        Then ${\cal X}\times_{{\mathbb C}^{n}} H$ can be described
\[
 (xy + \prod_{i=0}^{n}(z + s_i ) = 0)\wedge (\sum_{i=0}^{n}s_i = 0)
 \subset {\mathbb C}^{n+4}(x,y,z,s_0 , \cdots , s_n).
\]
        Let $u_0 := \sum_{i=0}^{n}s_i$,
	$u_i : = s_i - s_0$, $(1 \le i \le n)$
        and $z_i := z + s_i$. By using
        this coordinate, we can describe
\[
  {\cal X} \times_{{\mathbb C}^{n}} H = (xy + \prod_{i=0}^{n}z_i  = 0)
 \subset {\mathbb C}^{n+3}(x,y,z_0 , \cdots , z_n)
\]
        and $f = (z_1 - z_0 , \cdots , z_n - z_0 )$.
        Thus $f : {\cal X}\times_{{\mathbb C}^{n}}H 
               \to H$ is a Gorenstein homogeneous toric
        deformation.
 \item  Let $g : {\cal X} \to {\cal M}$ be a versal deformation space
        of $n$-dimensional $(n \ge 3)$
        isolated Gorenstein toric singularity. K.~Altmann
        shows in \cite[Theorem 8.1]{alt}
        that for every irreducible component ${\cal S}$ of ${\cal M}$,
\[
 \begin{CD}
   X= {\cal X}\times_{{\cal S}_{\text{red}}}{\cal M} @>>> {\cal X} \\
 @V{f}VV          @VV{g}V \\
   {\cal S}_{\text{red}} @>>>    {\cal M}     
 \end{CD}
\]
  the pull back 
  $f : X \to {\cal S}_{\text{red}}$ is a Gorenstein 
  homogeneous toric deformation.
\end{enumerate}

\section{Simultaneous terminalization for Gorenstein homogeneous
         toric deformation}
In this section, we give a proof of Theorem 1.
 Let $Y$ be a hypersurface in a quotient space such that
$$
 Y = (h := x_1 \cdots x_p - x_{p+1} \cdots x_{n+1} = 0) \subset 
 {\mathbb C}^{n+1}/G
$$
 where $G \cong {\mathbb Z}/n{\mathbb Z}$.
 The action of $G$ 
 on ${\mathbb C}^{n+1}$ is as follows:
$$
 (x_1 , \cdots , x_{n+1}) \to 
 (\zeta^{a_1}x_1 , \cdots , \zeta^{a_{n+1}}x_{n+1}), \quad (0 \le a_i < l)
$$
 where $\zeta$ is a $l$-th root of unity.
 We assume that 
 ${\mathbb C}^{n+1}/G$ has only Gorenstein terminal singularities.
 Moreover we assume that
$$
 \sum_{i=1}^{n+1} a_i \le l + \max \{\sum_{i=1}^{p} a_i , 
                                     \sum_{i=p+1}^{n+1} a_i \}.
$$
 Then $Y$ has canonical singularities by \cite[Theorem 4.6]{youngperson}.
 Let $X := (h = t) \subset {\mathbb C}^{n+1}/G \times {\mathbb C}(t)$
 and $f :X \to {\mathbb C}(t)$ projection. Then a general fibre
 of $f$ has only ${\mathbb Q}$-factorial terminal singularities
 because the total space 
 ${\mathbb C}^{n+1}/G \times {\mathbb C}$ has only ${\mathbb Q}$-factorial
 terminal singularities.
 Suppose
 that there exists a simultaneous terminalization after some
 finite base change. Let 
 $( h = t^m) \subset {\mathbb C}^{n+1}/G \times {\mathbb C}$ be
 a finite base change of $X$ and $ \tau : {\cal X} \to (h = t^m)$ 
 a simultaneous terminalization.
 We consider the following diagram:
$$
 \begin{array}{ccccc}
   (h = t) \subset {\mathbb C}^{n+2} 
  & \leftarrow & 
   (h = t^m)  \subset {\mathbb C}^{n+2}&
   \stackrel{\tau'}{\leftarrow} & {\cal X'} \\
   \downarrow & & \downarrow & & \downarrow \\
   X = (h = t) 
   \subset {\mathbb C}^{n+1}/G \times {\mathbb C} 
    &  \leftarrow & 
   (h = t^m)
   \subset {\mathbb C}^{n+1}/G \times {\mathbb C} 
    & \stackrel{\tau}{\leftarrow} & {\cal X}
 \end{array},
$$
 where ${\cal X}' = 
 {\cal X} \times_{{\mathbb C}^{n+1}/G \times {\mathbb C}} {\mathbb C}^{n+2}$.
 Because general fibres of $f : X \to {\mathbb C}$ have only 
 ${\mathbb Q}$-factorial terminal singularities, the codimension of
 exceptional sets of $\tau$ is greater than two. Thus $\tau$ and
 $\tau'$ are small birational morphisms. But 
 $(h =t^m )\subset {\mathbb C}^{n+2}$ has only  hypersurface singularities
 whose singular locus has codimension four. 
 Thus it is ${\mathbb Q}$-factorial,
 a contradiction. \qqed

\section{Simultaneous minimal model of Gorenstein homogeneous 
         toric deformation}
\begin{thm}
 Let $f : X := {\rm Spec}{\mathbb C}[\sigma^{\vee}\cap M] 
 \to {\mathbb C}^{m}$ be a Gorenstein homogeneous
 toric deformation and
 $\tau : \widetilde{X} \to X$ a toric minimal model of $X$.
 Assume that $\dim X = n+m$. Then
\begin{enumerate}
 \item $f \circ \tau : \widetilde{X} \to {\mathbb C}^{m}$ is a flat morphism.
 \item $K_{\widetilde{X}_t}$ is $\tau$-nef
 \item $\widetilde{X}_t$ has only canonical complete intersection
       singularities in quotient space such that 
\[
  (F_i - F_0 = 0) \subset 
  {\mathbb C}^{n+m}/G,  \quad (1 \le i \le m)
\]
 where $G$ is an abelian group acting on ${\mathbb C}^{n+m}$ 
 diagonally, ${\mathbb C}^{n+m}/G$ has only Gorenstein
 terminal singularities
 and $F_i$ are invariant monomials of the action of $G$.
 The monomials $F_i$ can be written as
\begin{eqnarray*}
 F_i &=& \prod_{j=p_i + 1}^{p_{i+1}} x_j \quad (0 \le i \le m) \\
     &&  0 = p_0  < p_1 < p_2 < \cdots < p_m < p_{m+1} = n+m
\end{eqnarray*}
 where $x_j$ is the j-th coordinate of  ${\mathbb C}^{n+m}$.
 \item If $\widetilde{X}$ is smooth, $\widetilde{X}_t$ has only
       terminal singularities.
\end{enumerate}
\end{thm}
{\sc Remark. \quad} If $\dim X = 2+m$ ( i.e. $f$ is 2-dimensional fibration ),
 then we can write $F_i$, $(1 \le i \le m)$ as
$$
 F_i = x_{i+1} \quad (0 \le i \le m-1) , \quad F_m = x_{m}x_{m+1}
$$
 by changing indices if necessary. Because $F_i$ are invariant monomials
 under the action of $G$, the action of each element of $g \in G$
 is nontrivial only $x_m$, $x_{m+1}$ axises. 
 But ${\mathbb C}^{2+m}/G$
 has only Gorenstein terminal singularity, the action of $G$ must be
 trivial. Thus each fibre of $f\circ \tau$ is smooth, and $\tau$ gives a
 simultaneous resolution.

%

\noindent
{\sc Proof of Theorem 3. \quad} 
 By K.~Altmann \cite[Theorem 3.5 , Remark 3.6]{homogeneous},
 we can state the construction of $\sigma$ as follows:
\begin{enumerate}
 \item $\sigma$ can be written as $\sigma = {\mathbb R}_{\ge 0}P$ where
 $P$ is an  $(n+m-1)$-dimensional polygon 
 such that
\begin{eqnarray*}
 P &:=& {\rm Conv}(\cup_{i=0}^{m}R_i \times e_i) \\
   &&\text{where $R_i \times e_i = 
       \{(x_1 , \cdots , x_{n-1}, 0 , \cdots , 1 , \cdots ,0)
          \in {\mathbb R}^{n+m}|
          (x_1 , \cdots , x_{n-1}) \in R_i \}$}
\end{eqnarray*}
 and $R_i$, $(0 \le i \le m)$ are integral polytopes in ${\mathbb R}^{n-1}$.
 \item
 $f$ can be written as $(x^{r_i} - x^{r_0})$, $(1 \le i \le m)$ where
 $r_i : N_{{\mathbb R}} = {\mathbb R}^{n+m} \to {\mathbb R}$ 
  is the $n+i$-th projection.
\end{enumerate}
 Thus, all primitive one dimensional generators of $\sigma$
 are contained in 
 the hyperplane defined by $r_0 + \cdots + r_m = 1$.
 Let $\tau : \widetilde{X} \to X$ be a toric minimal model of $X$ 
 and $\sigma = \cup \sigma_{\lambda}$ is the corresponding cone
 decomposition. Then these cones satisfy the following 
 three conditions:
\begin{enumerate}
 \item $\sigma_{\lambda}$ is a simplex.
 \item Let $ k_1 ,\cdots , k_{n+m} $ be
       one dimensional primitive generators of $\sigma_{\lambda}$.
       Then all $k_i$ are contained in the hypersurface 
       defined by $r_0 + \cdots + r_m = 1$.
 \item The polytope 
$$
 \Delta_{\lambda} := \sum_{i=0}^{n+m} \alpha_i k_i, \quad
 \sum_{i=0}^{n+m} \alpha_i \le 1, \quad \alpha_i \ge 0 
$$
       contains no lattice points except its vertices.
\end{enumerate}
 Let $X_{\lambda}:= {\rm Spec}{\mathbb C}[\sigma_{\lambda}^{\vee}\cap M]$
 and $k_{i}^{\vee}$, $(1 \le i \le n+m)$ be the dual vectors of $k_{i}$.
 Then, by (1), $X_{\lambda}$ can be written
 as follows:
$$
 X_{\lambda} \cong {\mathbb C}^{n+m}/G 
$$
 where 
 $G := N/\oplus_{i=1}^{n+m}{\mathbb Z}k_i$ and  the action of $G$ 
 is diagonal.
 Because $k_j$ are contained in the hypersurface defined by
 $r_0 + \cdots + r_m = 1$ and
 $\langle r_i , k_j  \rangle\ge 0$ ($r_i \in \sigma^{\vee}$), 
$$
\left\{
\begin{array}{cc}
  \langle r_i , k_j \rangle = 1 & \text{for $p_i   < j \le  p_{i+1}$} \\
  \langle r_i , k_j \rangle = 0 & \text{other $j$} 
\end{array}
\right.
$$
 where $0=p_0 < p_1 < p_2 < \cdots < p_{m} < p_{m+1} = n+m$. 
 Thus we can write
$$
 x^{r_i} = \prod_{j=p_{i}+1}^{p_{i+1}} x_{j}
$$
 where $x_j = x^{k^{\vee}_{j}}$ is the j-th coordinate of ${\mathbb C}^{n+m}$.
 The monomials $x^{r_i}$ are invariant under the action of $G$, 
 because $r_i \in \sigma_{\lambda}^{\vee}\cap M$. Thus
 if we set $F_i = x^{r_i}$,  proofs of theorem (1), (2) and 
 (3) are completed. Finally, we show (4). Since $X_{\lambda}$ is 
 smooth, a general fibre of $f\circ \tau$ is smooth. 
 We study a singularity of central fibre.
 Because $X_{\lambda}$ is smooth
 we may set $k_i = e^{i}$ $(1 \le i \le n+m)$ where
 $e^{i}$ are the standard basis of ${\mathbb Z}^{n+m}$. Then
 $r_i$ can be written 
$$
 r_i = (e^{p_i + 1})^{\vee} + \cdots + (e^{p_{i+1}})^{\vee}.
$$
 Thus a generator of the cone of central fibre $f\circ  \tau$
 can be written
$$
 e^{s_1} + \cdots + e^{s_m}, \quad p_i < s_i \le p_{i+1},
$$
 by the definition of homogeneous toric deformation.
 So all generators of this cone contains hypercube of ${\mathbb Z}^{n+m}$,
 hence central fibre has only terminal singularities.
\qqed

\noindent
A toric minimal model of $X$ is not unique, and which
can be joined by flops each other. We can obtain 
some examples of 4-dimensional flops as an  application.

\noindent
{\sc Example. \quad}
There are  4-dimensional flops which satisfy the following diagram:
$$
\begin{array}{ccc}
 {\mathbb P}(1,a,b) \subset Z & - \to & {\mathbb P}(1,a,b) 
  \subset Z^{+} \\
 \phi  \searrow & & \swarrow \phi^{+} \\
 & W &
\end{array}
$$
where $\phi({\mathbb P}(1,a,b)) = \phi^{+}({\mathbb P}(1,a,b)) = {\rm pt}$.

\noindent
{\sc Remark.} \quad In the case $a = b = 1$, it is known as ``Mukai
 transformation''.

\noindent
{\sc Construction of Example. \quad}
 Let $\bar{\sigma}$ be a 3-dimensional cone whose primitive generators
 are 
$$
 \bar{\sigma} = \langle (1,0,1),(0,1,1),(-a,-b,1),(1,1,1),(-a,1-b,1),
                        (1-a,-b,1) \rangle \subset {\mathbb R}^{3}.
$$
 We consider a deformation of toric singularity 
 $Y := {\rm Spec}{\mathbb C}[\bar{\sigma}^{\vee} \cap \bar{M}]$.
 The polytope $Q := \bar{\sigma}\cap 
                        \{(x,y,z) \in {\mathbb R}^{3} | z=1\}$ has
Minkowski decomposition $R_0 + R_1 + R_2$ where
$$
 R_0 = \langle (1,0),(0,0) \rangle, \quad
 R_1 = \langle (0,1),(0,0) \rangle, \quad
 R_2 = \langle (-a,-b),(0,0) \rangle.
$$
Thus by \cite[Remark 3.6]{homogeneous} the corresponding toric homogeneous 
deformation space of $Y$ is a toric variety
$X := {\rm Spec}{\mathbb C}[\sigma^{\vee} \cap M]$ where
$$
 \sigma := {\mathbb R}_{\ge 0}\langle e_1 , e_2 , e_3 , e_4 , e_5 , e_6 \rangle,
$$
$$
\begin{array}{ccc}
 e_1 :=  (1,0,1,0,0), & e_2 := (0,0,1,0,0), & e_3 := (0,1,0,1,0) \\
 e_4 :=  (0,0,0,1,0), & e_5 := (-a,-b,0,0,1), & e_6 := (0,0,0,0,1).
\end{array}
$$
 We construct two different crepant resolutions
 of $X$.
 Because $ae_1  -ae_2 + be_3 - be_4 + e_5 - e_6 = 0$, 
 by Reid\cite[Lemma 3.2]{reid}, there are 
 two cone decompositions of $\sigma$ such that
$$
 \sigma = \langle e_2 , e_4 , e_6, e_1 , e_3 \rangle \cup
          \langle e_2 , e_4 , e_6, e_3 , e_5 \rangle \cup
          \langle e_2 , e_4 , e_6, e_5 , e_1 \rangle
$$
 and 
$$ 
 \sigma =  \langle e_1 , e_3 , e_5 , e_2 , e_4 \rangle \cup
           \langle e_1 , e_3 , e_5 , e_4 , e_6 \rangle \cup
           \langle e_1 , e_3 , e_5 , e_6 , e_2 \rangle.
$$
 Let $\tilde{Z}$ and 
     $\tilde{Z^{+}}$ be the toric varieties corresponding 
     to above cone decompositions.
 Then the exceptional sets of $\phi : \tilde{Z} \to X$ and
 $\phi^{+} : \tilde{Z^{+}} \to X$ are isomorphic 
 to ${\mathbb P}(1,a,b)$ and $\phi({\mathbb P}(1,a,b)) 
 = \phi^{+}({\mathbb P}(1,a,b)) = pt$.
 There is a diagram
$$
 \begin{array}{ccccc}
  \tilde{Z} & \stackrel{\phi}{\to} & X  & 
              \stackrel{\phi^{+}}{\leftarrow} & \tilde{Z^{+}} \\
   f \downarrow & & \downarrow & & \downarrow f^{+}             \\
   {\mathbb C}^2 & = & {\mathbb C}^{2} & = & {\mathbb C}^{2}.
 \end{array}
$$
 Then the exceptional set of $\phi$ and $\phi^{+}$ are contained in a fibre
 of $f$ and $f^{+}$ respectively, since  $\phi({\mathbb P}(1,a,b)) 
 = \phi^{+}({\mathbb P}(1,a,b)) = pt$.
 Let $\imath : {\mathbb C} \to {\mathbb C}^{2}$ be the diagonal map.
 We set $Z$, $Z^{+}$ and 
 $W$  the pull back by $\imath$ of 
 $\tilde{Z} \to {\mathbb C}^{2}$,
 $\tilde{Z^{+}} \to {\mathbb C}^{2}$ and
 $X \to {\mathbb C}^{2}$ respectively.
\qqed

\end{document}